\documentclass{amsart}
\usepackage{amsmath}
\usepackage{amsfonts}

\newtheorem{theorem}{Theorem}
\theoremstyle{plain}

\newtheorem{corollary}{Corollary}

\newtheorem{definition}{Definition}

\newtheorem{lemma}{Lemma}

\newtheorem{remark}{Remark}

\numberwithin{equation}{section}
\input{tcilatex}

\begin{document}
\title[On new inequalities]{New generalization fractional inequalities of
Ostrowski-Gr\"{u}ss type }
\author{Mehmet Zeki Sarikaya$^{\star }$}
\address{Department of Mathematics, Faculty of Science and Arts, D\"{u}zce
University, D\"{u}zce, Turkey}
\email{sarikayamz@gmail.com}
\thanks{$^{\star }$corresponding author}
\author{Hatice YALDIZ}
\email{yaldizhatice@gmail.com}

\begin{abstract}
In this paper, we use the Riemann-Liouville fractional integrals to
establish some new integral inequalities of Ostrowski-Gr\"{u}ss type. From
our results, the classical Ostrowski-Gr\"{u}ss type inequalities can be
deduced as some special cases.
\end{abstract}

\maketitle

\section{Introduction}

Let $f:[a,b]\rightarrow \mathbb{R}$ be continuous on $\left[ a,b\right] $
and differentiable on $\left( a,b\right) $ and assume $\left\vert f^{\prime
}(x)\right\vert \leq M$ for all $x\in (a,b).$ Then the following inequality
holds:%
\begin{equation}
\left\vert S(f;a,b)\right\vert \leq \frac{M}{b-a}\left[ \left( \frac{b-a}{2}%
\right) ^{2}+\left( x-\frac{a+b}{2}\right) ^{2}\right]  \label{Z}
\end{equation}%
for all $x\in \left[ a,b\right] $ where%
\begin{equation*}
S(f;a,b)=f(x)-\mathcal{M}\left( f;a,b\right)
\end{equation*}%
and 
\begin{equation}
\mathcal{M}\left( f;a,b\right) =\frac{1}{b-a}\int_{a}^{b}f(x)dx.  \label{Z1}
\end{equation}%
This inequality is well known in the literature as Ostrowski inequality.

In 1882, P. L. \v{C}eby\v{s}ev \cite{Cebysev} gave the following inequality:%
\begin{equation}
\left\vert T(f,g)\right\vert \leq \frac{1}{12}(b-a)^{2}\left\Vert f^{\prime
}\right\Vert _{\infty }\left\Vert g^{\prime }\right\Vert _{\infty },
\label{Z2}
\end{equation}%
where $f,g:[a,b]\rightarrow \mathbb{R}$ are absolutely continuous function,
whose first derivatives $f^{\prime }$ and $g^{\prime }$ are bounded,%
\begin{eqnarray}
T(f,g) &=&\frac{1}{b-a}\int\limits_{a}^{b}f(x)g(x)dx-\left( \frac{1}{b-a}%
\int\limits_{a}^{b}f(x)dx\right) \left( \frac{1}{b-a}\int%
\limits_{a}^{b}g(x)dx\right)  \notag \\
&&  \label{Z3} \\
&=&\mathcal{M}\left( fg;a,b\right) -\mathcal{M}\left( f;a,b\right) \mathcal{M%
}\left( g;a,b\right)  \notag
\end{eqnarray}%
and $\left\Vert .\right\Vert _{\infty }$ denotes the norm in $L_{\infty
}[a,b]$ defined as $\left\Vert p\right\Vert _{\infty }=\underset{t\in
\lbrack a,b]}{ess\sup }\left\vert p(t)\right\vert .$

In 1935, G. Gr\"{u}ss \cite{Gruss} proved the following inequality:

\begin{equation}
\left\vert \frac{1}{b-a}\int\limits_{a}^{b}f(x)g(x)dx-\frac{1}{b-a}%
\int\limits_{a}^{b}f(x)dx\frac{1}{b-a}\int\limits_{a}^{b}g(x)dx\right\vert
\leq \frac{1}{4}(\Phi -\varphi )(\Gamma -\gamma ),  \label{Z4}
\end{equation}%
provided that $f$ and $g$ are two integrable function on $[a,b]$ satisfying
the condition%
\begin{equation}
\varphi \leq f(x)\leq \Phi \text{ \ and \ }\gamma \leq g(x)\leq \Gamma \text{
for all }x\in \lbrack a,b].  \label{Z5}
\end{equation}%
The constant $\frac{1}{4}$ is best possible.

From \cite{[13]}, if $f:[a,b]\rightarrow \mathbb{R}$ is differentiable on $%
[a,b]$ with the first derivative $f^{\prime }$ integrable on $[a,b],$ then
Montgomery identity holds:%
\begin{equation}
f(x)=\frac{1}{b-a}\int\limits_{a}^{b}f(t)dt+\int\limits_{a}^{b}P_{1}(x,t)f^{%
\prime }(t)dt,  \label{7}
\end{equation}%
where $P_{1}(x,t)$ is the Peano kernel defined by%
\begin{equation*}
P_{1}(x,t)=\left\{ 
\begin{array}{ll}
\dfrac{t-a}{b-a}, & a\leq t<x \\ 
\dfrac{t-b}{b-a}, & x\leq t\leq b.%
\end{array}%
\right.
\end{equation*}

This inequality provides an upper bound for the approximation of integral
mean of a function $f$ by the functional value $f(x)$ at $x\in \lbrack a,b].$
In 2001, Cheng \cite{[1]} proved the following Ostrowski-Gr\"{u}ss type
integral inequality.

\begin{theorem}
\label{t1} Let $I\subset 
\mathbb{R}
$ be an open interval, $a,b\in I,a<b$. If $f:I\rightarrow 
\mathbb{R}
$ is a differentiable function such that there exist constants $\gamma
,\Gamma \in 
\mathbb{R}
$, with $\varphi \leq f^{^{\prime }}\left( x\right) \leq \Phi $, $x\in \left[
a,b\right] $. Then have%
\begin{eqnarray}
&&\left\vert f\left( x\right) -\frac{f\left( b\right) -f\left( a\right) }{b-a%
}\left( x-\frac{a+b}{2}\right) -\frac{1}{b-a}\int\limits_{a}^{b}f\left(
t\right) dt\right\vert  \label{hh} \\
&&  \notag \\
&\leq &\frac{1}{4}\left( b-a\right) \left( \Phi -\varphi \right) \text{, for
all }x\in \left[ a,b\right] .  \notag
\end{eqnarray}
\end{theorem}

In \cite{[5]}, Matic et. al gave the following theorem by use of Gr\"{u}ss
inequality:

\begin{theorem}
\label{t2} Let the assumptions of Theorem \ref{t1} hold. Then for all $x\in %
\left[ a,b\right] $, we have
\end{theorem}

\begin{eqnarray}
&&\left\vert f\left( x\right) -\frac{f\left( b\right) -f\left( a\right) }{b-a%
}\left( x-\frac{a+b}{2}\right) -\frac{1}{b-a}\int\limits_{a}^{b}f\left(
t\right) dt\right\vert  \label{hhh} \\
&&  \notag \\
&\leq &\frac{1}{4\sqrt{3}}\left( b-a\right) \left( \Phi -\varphi \right) . 
\notag
\end{eqnarray}

In \cite{bernant}, Bernett et al., by the use of Chebeyshev's functional,
improved the Matic et al. result by providing first membership of the right
side of (\ref{hhh}) in terms of Euclidean norm as follows:

\begin{theorem}
Let $f:[a,b]\rightarrow 
\mathbb{R}
$ be an absolutely continuous mapping whose derivative $f^{\prime }\in
L_{2}[a,b]$, then we have,
\end{theorem}

\begin{eqnarray*}
&&\left\vert f\left( x\right) -\frac{f\left( b\right) -f\left( a\right) }{b-a%
}\left( x-\frac{a+b}{2}\right) -\frac{1}{b-a}\int\limits_{a}^{b}f\left(
t\right) dt\right\vert \\
&& \\
&\leq &\frac{\left( b-a\right) }{2\sqrt{3}}\left( \frac{1}{(b-a)}\left\Vert
f^{^{\prime }}\right\Vert _{2}^{2}-\left( \frac{f\left( b\right) -f\left(
a\right) }{(b-a)}\right) ^{2}\right) ^{\frac{1}{2}} \\
&& \\
&\leq &\frac{\left( b-a\right) \left( \Phi -\varphi \right) }{4\sqrt{3}},%
\text{ if }\varphi \leq f^{\prime }(x)\leq \Phi \text{ for a.e.t on }[a,b],
\end{eqnarray*}%
for all $x\in \left[ a,b\right] .$

During the past few years many researchers have given considerable attention
to the above inequalities and various generalizations, extensions and
variants of these inequalities have appeared in the literature, see \cite%
{[1]}-\cite{bernant}, \cite{[5]}, \cite{[24]} and the references cited
therein. For recent results and generalizations concerning Ostrowski and Gr%
\"{u}ss inequalities, we refer the reader to the recent papers \cite{[1]}-%
\cite{bernant}, \cite{[5]}, \cite{zafar}-\cite{[24]}.

The theory of fractional calculus has known an intensive development over
the last few decades. It is shown that derivatives and\ integrals of
fractional type provide an adequate mathematical modelling of real objects\
and processes see (\cite{[6]}-\cite{sarikaya}). Therefore, the study of
fractional differential equations need more developmental of inequalities of
fractional type. The main aim of this work is to develop new integral
inequalities of Ostrowski-Gr\"{u}ss type for Riemann-Liouville fractional
integrals. From our results, the classical Ostrowski-Gr\"{u}ss type
inequalities can be deduced as some special cases. Let us begin by
introducing this type of inequality.

In \cite{[6]} and \cite{[21]}, the authors established some inequalities for
differentiable mappings which are connected with Ostrowski type inequality
by used the Riemann-Liouville fractional integrals, and they used the
following lemma to prove their results:

\begin{lemma}
\label{l} Let $f:I\subset \mathbb{R}\rightarrow \mathbb{R}$ be
differentiable function on $I^{\circ }$ with $a,b\in I$ ($a<b$) and $%
f^{\prime }\in L_{1}[a,b]$, then%
\begin{equation}
f(x)=\frac{\Gamma (\alpha )}{b-a}(b-x)^{1-\alpha }{\Large J}_{a}^{\alpha
}f(b)-{\Large J}_{a}^{\alpha -1}(P_{2}(x,b)f(b))+{\Large J}_{a}^{\alpha
}(P_{2}(x,b)f^{^{\prime }}(b)),\ \ \ \alpha \geq 1,  \label{z}
\end{equation}%
where $P_{2}(x,t)$ is the fractional Peano kernel defined by%
\begin{equation}
P_{2}(x,t)=\left\{ 
\begin{array}{ll}
\left( \dfrac{t-a}{b-a}\right) (b-x)^{1-\alpha }\Gamma (\alpha ), & a\leq t<x
\\ 
&  \\ 
\left( \dfrac{t-b}{b-a}\right) (b-x)^{1-\alpha }\Gamma (\alpha ), & x\leq
t\leq b.%
\end{array}%
\right.  \label{z1}
\end{equation}
\end{lemma}

In \cite{[6]} and \cite{[21]}, the authors derive the following interesting
fractional integral inequality:

\begin{eqnarray*}
&&\left\vert f\left( x\right) -\frac{1}{b-a}\left( b-x\right) ^{1-\alpha
}\Gamma \left( \alpha \right) J_{a}^{\alpha }\left( f\left( b\right) \right)
+J_{a}^{\alpha -1}\left( P_{2}\left( x,b\right) f\left( b\right) \right)
\right\vert \\
&& \\
&\leq &\frac{M}{\alpha \left( \alpha +1\right) }\left[ \left( b-x\right)
\left( 2\alpha \left( \frac{b-x}{b-a}\right) -\alpha -1\right) +\left(
b-a\right) ^{\alpha }\left( b-x\right) ^{1-\alpha }\right]
\end{eqnarray*}%
under the assumption that $\left\vert f^{\prime }\left( x\right) \right\vert
\leq M$, for any $x\in \left[ a,b\right] $.

Firstly, we give some necessary definitions and mathematical preliminaries
of fractional calculus theory which are used further in this paper. More
details, one can consult \cite{[17]} and \cite{[12]}.

\begin{definition}
The Riemann-Liouville fractional integral operator of order $\alpha \geq 0$
with $a\geq 0$ is defined as%
\begin{eqnarray*}
J_{a}^{\alpha }f(x) &=&\frac{1}{\Gamma (\alpha )}\dint\limits_{a}^{x}(x-t)^{%
\alpha -1}f(t)dt, \\
J_{a}^{0}f(x) &=&f(x).
\end{eqnarray*}
\end{definition}

Recently, many authors have studied a number of inequalities by used the
Riemann-Liouville fractional integrals, see (\cite{[6]}-\cite{sarikaya}) and
the references cited therein.

\section{Main Results}

\begin{theorem}
\label{t} Let $f:I\subseteq 
\mathbb{R}
\rightarrow 
\mathbb{R}
$ be a differentiable mapping in $I^{0}$ (interior of $I$), and $a,b\in
I^{0} $ with $a<b$ and $f^{\prime }\in L_{2}[a,b]$. If $f^{^{\prime
}}:\left( a,b\right) \rightarrow 
\mathbb{R}
$ is bounded on $\left( a,b\right) $ with $\varphi \leq f^{\prime }(x)\leq
\Phi $, then we have%
\begin{eqnarray}
&&\left\vert \frac{f\left( x\right) }{\Gamma \left( \alpha \right) }-\frac{%
\left( b-x\right) ^{1-\alpha }}{\left( b-a\right) }J_{a}^{\alpha }f\left(
b\right) +\frac{1}{\Gamma \left( \alpha \right) }J_{a}^{\alpha -1}\left(
P\left( x,b\right) f\left( b\right) \right) \right.  \notag \\
&&  \notag \\
&&\left. -\left( \frac{f\left( b\right) -f\left( a\right) }{b-a}\right)
\times \left( \frac{\left( b-x\right) ^{1-\alpha }\left( b-a\right) ^{\alpha
}}{\Gamma \left( \alpha +2\right) }-\frac{\left( b-x\right) }{\Gamma \left(
\alpha +1\right) }\right) \right\vert  \notag \\
&&  \label{h} \\
&\leq &(b-a)\left( K(x)\right) ^{\frac{1}{2}}\left( \frac{1}{(b-a)\Gamma
^{2}\left( \alpha \right) }\left\Vert f^{^{\prime }}\right\Vert
_{2}^{2}-\left( \frac{f\left( b\right) -f\left( a\right) }{(b-a)\Gamma
\left( \alpha \right) }\right) ^{2}\right) ^{\frac{1}{2}}  \notag \\
&&  \notag \\
&\leq &\frac{\left( K(x)\right) ^{\frac{1}{2}}}{2\Gamma \left( \alpha
\right) }(b-a)(\Phi -\varphi )  \notag
\end{eqnarray}%
for all $x\in \lbrack a,b]$ and $\alpha \geq 1$ where%
\begin{eqnarray*}
K(x) &=&\left( b-x\right) ^{1-\alpha }\left( b-a\right) ^{2\alpha -2}\left( 
\frac{1}{2\alpha +1}+\frac{1}{2\alpha -1}-\frac{1}{\alpha }\right) \\
&& \\
&&+\frac{(b-x)^{\alpha }}{\left( b-a\right) ^{2}}\left( \frac{b-x}{\alpha }-%
\frac{b-a}{2\alpha -1}\right) -\left( \frac{\left( b-x\right) ^{1-\alpha
}\left( b-a\right) ^{\alpha -1}}{\alpha \left( \alpha +1\right) }-\frac{%
\left( b-x\right) }{\alpha \left( b-a\right) }\right) ^{2}.
\end{eqnarray*}
\end{theorem}

\begin{proof}
We consider the fractional Peano kernel $P_{2}:\left[ a,b\right]
^{2}\rightarrow 
\mathbb{R}
$ as defined in (\ref{z1}). Using Korkine's identity%
\begin{equation*}
T\left( f,g\right) :=\frac{1}{2\left( b-a\right) ^{2}}\dint\limits_{a}^{b}%
\dint\limits_{a}^{b}\left( f\left( t\right) -f\left( s\right) \right) \left(
g\left( t\right) -g\left( s\right) \right) dsdt
\end{equation*}%
we obtain%
\begin{eqnarray*}
&&\frac{1}{\left( b-a\right) \Gamma ^{2}(\alpha )}\dint\limits_{a}^{b}\left(
b-t\right) ^{\alpha -1}P_{2}\left( x,t\right) f^{^{\prime }}\left( t\right)
dt \\
&&-\left( \frac{1}{\left( b-a\right) \Gamma (\alpha )}\dint\limits_{a}^{b}%
\left( b-t\right) ^{\alpha -1}P_{2}\left( x,t\right) dt\right) \left( \frac{1%
}{\left( b-a\right) \Gamma (\alpha )}\dint\limits_{a}^{b}f^{^{\prime
}}\left( t\right) dt\right) \\
&& \\
&=&\frac{1}{2\left( b-a\right) ^{2}\Gamma ^{2}(\alpha )}\dint\limits_{a}^{b}%
\dint\limits_{a}^{b}\left[ \left( b-t\right) ^{\alpha -1}P_{2}\left(
x,t\right) -\left( b-s\right) ^{\alpha -1}P_{2}\left( x,s\right) \right] %
\left[ f^{^{\prime }}\left( t\right) -f^{^{\prime }}\left( s\right) \right]
dsdt
\end{eqnarray*}%
i.e.%
\begin{eqnarray}
&&\frac{1}{\left( b-a\right) \Gamma \left( \alpha \right) }J_{a}^{\alpha
}\left( P_{2}\left( x,b\right) f^{^{\prime }}\left( b\right) \right) -\frac{1%
}{\left( b-a\right) \Gamma \left( \alpha \right) }J_{a}^{\alpha }\left(
P_{2}\left( x,b\right) \right) \left( \frac{f\left( b\right) -f\left(
a\right) }{b-a}\right)  \notag \\
&&  \label{h1} \\
&=&\frac{1}{2\left( b-a\right) ^{2}\Gamma ^{2}(\alpha )}\dint\limits_{a}^{b}%
\dint\limits_{a}^{b}\left[ \left( b-t\right) ^{\alpha -1}P_{2}\left(
x,t\right) -\left( b-s\right) ^{\alpha -1}P_{2}\left( x,s\right) \right] %
\left[ f^{^{\prime }}\left( t\right) -f^{^{\prime }}\left( s\right) \right]
dsdt.  \notag
\end{eqnarray}%
Since%
\begin{equation}
J_{a}^{\alpha }\left( P_{2}\left( x,b\right) f^{^{\prime }}\left( b\right)
\right) =f\left( x\right) -\frac{\Gamma \left( \alpha \right) }{b-a}\left(
b-x\right) ^{1-\alpha }J_{a}^{\alpha }f\left( b\right) +J_{a}^{\alpha
-1}\left( P_{2}\left( x,b\right) f\left( b\right) \right)  \label{h2}
\end{equation}%
and%
\begin{equation}
J_{a}^{\alpha }\left( P_{2}\left( x,b\right) \right) =\frac{\left(
b-x\right) ^{1-\alpha }\left( b-a\right) ^{\alpha }}{\alpha \left( \alpha
+1\right) }-\frac{\left( b-x\right) }{\alpha },  \label{h3}
\end{equation}%
then by (\ref{h1}) we get the following identity,%
\begin{eqnarray}
&&\frac{f\left( x\right) }{\left( b-a\right) \Gamma \left( \alpha \right) }-%
\frac{\left( b-x\right) ^{1-\alpha }}{\left( b-a\right) ^{2}}J_{a}^{\alpha
}f\left( b\right) +\frac{1}{\left( b-a\right) \Gamma \left( \alpha \right) }%
J_{a}^{\alpha -1}\left( P\left( x,b\right) f\left( b\right) \right)  \notag
\\
&&  \notag \\
&&-\frac{1}{(b-a)}\left( \frac{f\left( b\right) -f\left( a\right) }{b-a}%
\right) \times \left( \frac{\left( b-x\right) ^{1-\alpha }\left( b-a\right)
^{\alpha }}{\Gamma \left( \alpha +2\right) }-\frac{\left( b-x\right) }{%
\Gamma \left( \alpha +1\right) }\right)  \notag \\
&&  \label{h4} \\
&=&\frac{1}{2\left( b-a\right) ^{2}\Gamma ^{2}(\alpha )}\dint\limits_{a}^{b}%
\dint\limits_{a}^{b}\left[ \left( b-t\right) ^{\alpha -1}P_{2}\left(
x,t\right) -\left( b-s\right) ^{\alpha -1}P_{2}\left( x,s\right) \right] %
\left[ f^{^{\prime }}\left( t\right) -f^{^{\prime }}\left( s\right) \right]
dsdt.  \notag
\end{eqnarray}%
Using the Cauchy-Swartz inequality for double integrals, we write%
\begin{eqnarray}
&&\frac{1}{2\left( b-a\right) ^{2}\Gamma ^{2}(\alpha )}\left\vert
\dint\limits_{a}^{b}\dint\limits_{a}^{b}\left[ \left( b-t\right) ^{\alpha
-1}P_{2}\left( x,t\right) -\left( b-s\right) ^{\alpha -1}P_{2}\left(
x,s\right) \right] \left[ f^{^{\prime }}\left( t\right) -f^{^{\prime
}}\left( s\right) \right] dsdt\right\vert  \notag \\
&&  \notag \\
&\leq &\left( \frac{1}{2\left( b-a\right) ^{2}\Gamma ^{2}(\alpha )}%
\dint\limits_{a}^{b}\dint\limits_{a}^{b}\left[ \left( b-t\right) ^{\alpha
-1}P_{2}\left( x,t\right) -\left( b-s\right) ^{\alpha -1}P_{2}\left(
x,s\right) \right] ^{2}dsdt\right) ^{\frac{1}{2}}  \notag \\
&&  \label{h5} \\
&&\times \left( \frac{1}{2\left( b-a\right) ^{2}\Gamma ^{2}(\alpha )}%
\dint\limits_{a}^{b}\dint\limits_{a}^{b}\left[ f^{^{\prime }}\left( t\right)
-f^{^{\prime }}\left( s\right) \right] ^{2}dsdt\right) ^{\frac{1}{2}}. 
\notag
\end{eqnarray}%
However,%
\begin{eqnarray}
&&\frac{1}{2\left( b-a\right) ^{2}\Gamma ^{2}(\alpha )}\dint\limits_{a}^{b}%
\dint\limits_{a}^{b}\left[ \left( b-t\right) ^{\alpha -1}P_{2}\left(
x,t\right) -\left( b-s\right) ^{\alpha -1}P_{2}\left( x,s\right) \right]
^{2}dsdt  \notag \\
&&  \notag \\
&=&\frac{1}{\left( b-a\right) \Gamma ^{2}(\alpha )}\dint\limits_{a}^{b}%
\left( b-t\right) ^{2\alpha -2}P_{2}^{2}\left( x,t\right) dt-\left( \frac{1}{%
\left( b-a\right) \Gamma (\alpha )}\dint\limits_{a}^{b}\left( b-t\right)
^{\alpha -1}P_{2}\left( x,t\right) dt\right) ^{2}  \notag \\
&&  \label{h6} \\
&=&\left( b-x\right) ^{1-\alpha }\left( b-a\right) ^{2\alpha -2}\left( \frac{%
1}{2\alpha +1}+\frac{1}{2\alpha -1}-\frac{1}{\alpha }\right) +\frac{\left(
b-x\right) ^{\alpha }}{\left( b-a\right) ^{2}}\left( \frac{b-x}{\alpha }-%
\frac{b-a}{2\alpha -1}\right)  \notag \\
&&  \notag \\
&&-\left( \frac{\left( b-x\right) ^{1-\alpha }\left( b-a\right) ^{\alpha -1}%
}{\alpha \left( \alpha +1\right) }-\frac{\left( b-x\right) }{\alpha \left(
b-a\right) }\right) ^{2},  \notag
\end{eqnarray}%
and%
\begin{equation}
\frac{1}{2\left( b-a\right) ^{2}\Gamma ^{2}\left( \alpha \right) }%
\dint\limits_{a}^{b}\dint\limits_{a}^{b}\left( f^{^{\prime }}\left( t\right)
-f^{^{\prime }}\left( s\right) \right) ^{2}dsdt=\frac{1}{(b-a)\Gamma
^{2}\left( \alpha \right) }\left\Vert f^{^{\prime }}\right\Vert
_{2}^{2}-\left( \frac{f\left( b\right) -f\left( a\right) }{(b-a)\Gamma
\left( \alpha \right) }\right) ^{2}.  \label{h7}
\end{equation}%
Using (\ref{h5})-(\ref{h7}), we deduce the (\ref{h1}) inequality.

Moreover, if $\varphi \leq f^{\prime }(t)\leq \Phi $ almost everywhere $t$
on $(a,b),$ then by using Gr\"{u}ss inequality, we get%
\begin{equation*}
0\leq \frac{1}{b-a}\dint\limits_{a}^{b}\left( f^{\prime }\left( t\right)
\right) ^{2}dt-\left( \frac{1}{b-a}\dint\limits_{a}^{b}f^{\prime }\left(
t\right) dt\right) ^{2}\leq \frac{1}{4}(\Phi -\varphi )^{2},
\end{equation*}%
which proves the last inequality of (\ref{h}).
\end{proof}

\begin{corollary}
\label{c} Under the assumptions of Theorem \ref{t} with $\alpha =1$. Then
the following inequality holds:For 
\begin{eqnarray}
&&\left\vert f\left( x\right) -\frac{1}{b-a}\dint\limits_{a}^{b}f\left(
t\right) dt-\frac{f\left( b\right) -f\left( a\right) }{b-a}\left( x-\frac{a+b%
}{2}\right) \right\vert  \label{c1} \\
&&  \notag \\
&\leq &\frac{(b-a)}{2\sqrt{3}}\left( \frac{1}{(b-a)}\left\Vert f^{^{\prime
}}\right\Vert _{2}^{2}-\left( \frac{f\left( b\right) -f\left( a\right) }{%
(b-a)}\right) ^{2}\right) ^{\frac{1}{2}}  \notag \\
&&  \notag \\
&\leq &\frac{(b-a)(\Phi -\varphi )}{4\sqrt{3}}.  \notag
\end{eqnarray}
\end{corollary}

\begin{proof}
Proof of Corollary \ref{c} can be as similar to the proof of Theorem \ref{t}.
\end{proof}

\begin{remark}
If we take $x=\frac{a+b}{2}$ in (\ref{c1}), it follows that 
\begin{eqnarray*}
&&\left\vert f\left( \frac{a+b}{2}\right) -\frac{1}{b-a}\dint%
\limits_{a}^{b}f\left( t\right) dt\right\vert \leq \frac{(b-a)}{2\sqrt{3}}%
\left( \frac{1}{(b-a)}\left\Vert f^{^{\prime }}\right\Vert _{2}^{2}-\left( 
\frac{f\left( b\right) -f\left( a\right) }{(b-a)}\right) ^{2}\right) ^{\frac{%
1}{2}} \\
&& \\
&\leq &\frac{(b-a)(\Phi -\varphi )}{4\sqrt{3}}.
\end{eqnarray*}
\end{remark}


\begin{thebibliography}{99}
\bibitem{[1]} X. L. Cheng, \textit{Improvement of some Ostrowski-Gr\"{u}ss
type inequalities}, Computers Math. Applic, 42 (2001), 109-114.

\bibitem{[2]} S. S. Dragomir and S. Wang,\textit{\ An inequality of
Ostrowski-Gr\"{u}ss type and its applications to the estimation of error
bounds for some special means and for home numerical quadrature rules},
Computers Math. Applic,33(11)(1997),15-20.

\bibitem{bernant} Barnett N. S., Dragomir S. S. and Sofo A., \textit{Better
bounds for an inequality of Ostrowski type with applications,} Demonstratio
Math., 34(3), 2001, 533-542.

\bibitem{Cebysev} P. L. \v{C}eby\v{s}ev, \textit{Sur less expressions
approximatives des integrales definies par les autres prises entre les memes
limites, }Proc. Math. Soc. Charkov, 2, 93-98, 1882.

\bibitem{Gruss} G. Gr\"{u}ss, \textit{\"{U}ber das maximum des absoluten
Betrages von }$\tfrac{1}{b-a}\int\limits_{a}^{b}{\small f(x)g(x)dx-}\frac{1}{%
(b-a)^{2}}\int\limits_{a}^{b}{\small f(x)dx}\int\limits_{a}^{b}{\small g(x)dx%
},$ Math. Z., 39, 215-226, 1935.

\bibitem{[5]} M. Matic, J. Pecaric and N. Ujevic, \textit{Improvement and
further generalization of inequalities of Ostrowski-Gr\"{u}ss type,}%
Computers Math. Appl., 39(3/4), (2000), 161-175.

\bibitem{[6]} G. Anastassiou, M.R. Hooshmandasl, A. Ghasemi and F.
Moftakharzadeh, \textit{Montgomery identities for fractional integrals and
related fractional inequalities}, J. Inequal. in Pure and Appl. Math, 10(4),
2009, Art. 97, 6 pp.

\bibitem{[7]} S. Belarbi and Z. Dahmani, \textit{On some new fractional
integral inequalities}, J. Inequal. in Pure and Appl. Math, 10(3), 2009,
Art. 97, 6 pp.

\bibitem{[9]} Z. Dahmani, L. Tabharit and S. Taf, \textit{Some fractional
integral inequalities}, Nonlinear Science Letters A, 2(1), 2010, p.155-160.

\bibitem{[10]} Z. Dahmani, L. Tabharit and S. Taf, \textit{New inequalities
via Riemann-Liouville fractional integration}, J. Advance Research Sci.
Comput., 2(1), 2010, p.40-45.

\bibitem{[21]} M.Z. Sarikaya and H. Ogunmez, \textit{On new inequalities via
Riemann-Liouville fractional integration}, arXiv:1005.1167v1, submitted.

\bibitem{[22]} M.Z. Sarikaya, E. Set, H. Yaldiz and N., Basak, \textit{%
Hermite -Hadamard's inequalities for fractional integrals and related
fractional inequalities}, Mathematical and Computer Modelling,
DOI:10.1016/j.mcm.2011.12.048.

\bibitem{sarikaya} M.Z. Sarikaya, H. Yaldiz and N. Basak, \textit{New
fractional inequalities of Ostrowski-Gr\"{u}ss type}, submitted.

\bibitem{[17]} S. G. Samko, A. A Kilbas, O. I. Marichev, \textit{\
Fractional Integrals and Derivatives Theory and Application}, Gordan and
Breach Science, New York, 1993.

\bibitem{[12]} R. Gorenflo and F. Mainardi, \textit{Fractionalcalculus:
integral and differentiable equations of fractional order}, Springer Verlag,
Wien, 1997, p.223-276.

\bibitem{[13]} D. S. Mitrinovic, J. E. Pecaric and A. M. Fink, \textit{%
Inequalities involving functions and their integrals and derivatives,}
Kluwer Academic Publishers, Dordrecht, 1991.

\bibitem{zafar} F. Zafar and N.A. Mir, \textit{A generalization of
Ostrowski-Gr\"{u}ss type inequality for first diffrentiable mappings,}
Tamsui Oxford J. of Math. Sci. 26(1), 2010, 61-76.

\bibitem{[23]} A. Rafiq and F. Ahmad, \textit{Another weighted Ostrowski-Gr%
\"{u}ss type inequality for twice differentiable mappings}, Kragujevac
Journal of Mathematics, 31 (2008), 43-51.

\bibitem{[24]} F. Tong and L. Guan, \textit{A simple proof of the
generalized Ostrowski-Gr\"{u}ss type integral inequality,} Int. Jour. of
Math. Analysis, 2(18), (2008), 889-892.
\end{thebibliography}
\end{document}